\documentclass[10pt]{amsart}

\usepackage{amsmath, amscd}
\usepackage{amsfonts}
\usepackage{amssymb}
 \usepackage{amsthm}
 \newtheorem{definition}{Definition}[section]
\newtheorem{theorem}{Theorem}[section]
\newtheorem{lemma}{Lemma}[section]

\newtheorem{example}{Example}[section]

\newcommand{\st}{^*}

\newcommand{\sh}{^{\sharp}}

\begin{document}

%\begin{frontmatter}

%% Title, authors and addresses

%% use the tnoteref command within \title for footnotes;
%% use the tnotetext command for the associated footnote;
%% use the fnref command within \author or \address for footnotes;
%% use the fntext command for the associated footnote;
%% use the corref command within \author for corresponding author footnotes;
%% use the cortext command for the associated footnote;
%% use the ead command for the email address,
%% and the form \ead[url] for the home page:
%%
%% \title{Title\tnoteref{label1}}
%% \tnotetext[label1]{}
%% \author{Name\corref{cor1}\fnref{label2}}
%% \ead{email address}
%% \ead[url]{home page}
%% \fntext[label2]{}
%% \cortext[cor1]{}
%% \address{Address\fnref{label3}}
%% \fntext[label3]{}

\title{Poisson Reduction}

%% use optional labels to link authors explicitly to addresses:
%% \author[label1,label2]{<author name>}
%% \address[label1]{<address>}
%% \address[label2]{<address>}

\author{Chiara Esposito}

\address{Department of Mathematics, Universitat Autonoma de Barcelona, 08193 Bellaterra. Spain}
\maketitle

\begin{abstract}
In this paper we develope a theory of reduction for classical systems with Poisson Lie groups symmetries using the notion of momentum map introduced by Lu. The local description of Poisson manifolds and Poisson Lie groups and the properties of Lu's momentum map allow us to define a Poisson reduced space.
\end{abstract}

%\begin{keyword}
%Momentum map \sep Poisson manifold \sep Poisson Lie group \sep group action \sep reduction
%% keywords here, in the form: keyword \sep keyword

%% MSC codes here, in the form: \MSC code \sep code
%% or \MSC[2008] code \sep code (2000 is the default)

%\end{keyword}

%\end{frontmatter}

%%
%% Start line numbering here if you want
%%
% \linenumbers

%% main text
\section{Introduction}\label{intro}
In this paper we prove a generalization of the Marsden-Weinstein reduction to the general case of an arbitrary Poisson Lie group action on a Poisson manifold. Reduction procedures are known in many different settings. In particular, a reduction theory is known in the case of Poisson Lie groups acting on symplectic manifolds \cite{Lu3} and in the case of Lie groups acting on Poisson manifolds \cite{RO}, \cite{MR}. An important generalization to the Dirac setting has been studied in \cite{BC}.

The theory of symplectic reduction plays a key role in classical mechanics. The phase 
space of a system of $n$ particles is described by a symplectic or more generally Poisson manifold. Given a symmetry group of dimension $k$ acting on a mechanical system, the dimension of the phase space can be reduced by $2k$. Marsden-Weinstein reduction formalizes this feature. Recall roughly the notion of Hamiltonian actions in this setting. Given a Poisson manifold $M$ there are natural Hamiltonian vector fields $\{f, \cdot \}$ on $M$. Let $G$ be a Lie group acting on $M$ by $\Phi$; the action is Hamiltonian if the vector fields defined by the infinitesimal generator of $\Phi$ are Hamiltonian. More precisely, let $G$ be a Lie group acting on a Poisson manifold $(M, \pi)$. The action $\Phi:G\times M \to M$ is canonical if it preserves the Poisson structure $\pi$. Suppose that there exists a linear map $H: \mathfrak{g} \to C^{\infty}(M)$ such that the infinitesimal generator $\Phi_{X}$ for $X\in \mathfrak{g}$ of the canonical action is induced by $H$ by
$$
\Phi_{X}= \{H_{X}, \cdot\}.
$$
A canonical action induced by $H$ is said Hamiltonian if $H$ is a Lie algebra homomorphism. We can define a map $\boldsymbol{\mu}: M \to \mathfrak{g}^*$, called momentum map, by $H_{X}(m)=\langle \boldsymbol{\mu}(m), X\rangle$ for $m\in M$. It is equivariant if the corresponding $H$ is a Lie algebra homomorphism. Given an Hamiltonian action, under certain assumptions, the reduced space has been defined as $M//G:=\boldsymbol{\mu}^{-1}(u)/G_{u}$ and it has been proved that it is a Poisson manifold \cite{MsWe}.

In this paper we are interested in analyzing the case in which one has an extra structure on the Lie group, a Poisson structure making it a Poisson Lie group. Poisson Lie groups are very interesting objects in mathematical physics. They may be regarded as classical limit of quantum groups \cite{Dr1} and they have been studied as carrier spaces of dynamical systems \cite{LMS}. It is believed that actions of Poisson Lie groups on Poisson manifolds should be used to understand the ``hidden symmetries'' of certain integrable systems \cite{STS}. Moreover, the study of classical systems with Poisson Lie group symmetries may give information about the corresponding quantum group invariant system (an attempt can be found in \cite{me}, \cite{me1}).

The purpose of this paper is to prove that, given a Poisson manifold acted by a Poisson Lie group, under certain conditions, we can also reduce this phase space to another Poisson manifold.

The paper is organized as follows. In Section \ref{sec_pg} we recall some basic elements of Poisson geometry: Poisson manifolds and their local description, Lie bialgebras and Poisson Lie groups. A nice review of these results can be found in \cite{V} and \cite{YK}. The Section \ref{sec_mm} is devoted to Poisson actions and associated momentum maps and we discuss dressing actions and their properties. In Section \ref{sec: pr} we present the main result of this paper, the Poisson reduction, and we discuss an example.

\noindent{\bf Acknowledgments:} I would like to thank my advisor Ryszard Nest and Eva Miranda for many interesting discussions about Poisson reduction and its possible developments. I also wish to thank George M. Napolitano for his help and his useful suggestions and Rui L. Fernandes for his comments regarding Dirac reduction theory.

\section{Poisson manifolds, Poisson Lie groups and Lie bialgebras}\label{sec_pg}

In this section we introduce the notion of Poisson manifolds and their local description, we give some background about Poisson Lie groups and Lie bialgebras which will be used in the paper. For more details on this subject, see \cite{Lu3}, \cite{Dr1}, \cite{YK}, \cite{V}, \cite{We1}.

\subsection{Poisson manifolds and symplectic foliation}\label{sec_1.1}

A Poisson structure on a smooth manifold $M$ is a Lie bracket $\{\cdot, \cdot\}$ on the space $C^{\infty}(M)$ of smooth functions on $M$ which satisfies the Leibniz rule. This bracket is called Poisson bracket and a manifold $M$ equipped with such a bracket is called Poisson manifold.
Therefore, a bivector field $\pi$ on $M$ such that the bracket
$$
\{ f, g\}:= \langle \pi, df\wedge dg\rangle
$$
is a Poisson bracket is called Poisson tensor or Poisson bivector field. A Poisson tensor can be regarded as a bundle map $\pi^{\sharp}: T^*M\to TM$:
$$
 \langle \alpha, \pi^{\sharp}(\beta)\rangle = \pi(\alpha,\beta)
$$

\begin{definition}
A mapping $\phi: (M_1,\pi_1)\rightarrow (M_2,\pi_2)$ between two Poisson manifolds is called a Poisson mapping if $\forall f,g\in C^{\infty}(M_2)$ one has
\begin{equation}
\lbrace f\circ\phi, g\circ\phi\rbrace_1=\lbrace f,g\rbrace_2\circ \phi
\end{equation}
\end{definition}
The structure of a Poisson manifold is described by the splitting theorem of Alan Weinstein \cite{We1}, which shows that locally a Poisson manifold is a direct product of a symplectic manifold with another Poisson manifold whose Poisson tensor vanishes at a point. 

\begin{theorem}[Weinstein]\label{thm: split}
On a Poisson manifold $(M,\pi)$, any point $m\in M$ has a coordinate neighborhood with coordinates
$(q_1,\dots,q_k,p_1,\dots,p_k,\allowbreak y_1,\dots,y_l)$ centered at $m$, such that
\begin{equation}\label{eq: splitp}
\pi=\sum_i \frac{\partial}{\partial q_i}\wedge\frac{\partial}{\partial p_i}+\frac{1}{2}\sum_{i,j}\phi_{ij}(y) \frac{\partial}{\partial y_i}\wedge\frac{\partial}{\partial y_j}\qquad \phi_{ij}(0)=0.
\end{equation}
The rank of $\pi$ at $m$ is $2k$. Since $\phi$ depends only on the $y_i$s, this theorem gives
a decomposition of the neighborhood of $m$ as a product of two Poisson manifolds: one with rank $2k$, and the other with rank 0 at $m$.
\end{theorem}
 The term
 \begin{equation}
\frac{1}{2}\sum_{i,j}\phi_{ij}(y)\frac{\partial}{\partial y_i}\wedge \frac{\partial}{\partial y_j}
\end{equation}
is called transverse Poisson structure and it is evident that the equations $y_{i}=0$ determine the symplectic leaf through $m$.

\subsection{Lie bialgebras and Poisson Lie groups}
\label{sec:lie bialgebras}

\begin{definition}
A Poisson Lie group $(G,\pi_G)$ is a Lie group equipped with a multiplicative Poisson structure $\pi_G$, i.e. such that the multiplication map $G\times G \to G$ is a Poisson map.
\end{definition}

Let $G$ be a Lie group with Lie algebra $\mathfrak{g}$. The linearization  $\delta:= d_{e}\pi_{G}: \mathfrak{g}\to \mathfrak{g}\wedge \mathfrak{g}$ of $\pi_{G}$ at $e$ defines a Lie algebra structure on the dual $\mathfrak{g}^*$ of $\mathfrak{g}$ and, for this reason, it is called cobracket. The pair $(\mathfrak{g},\mathfrak{g}^*)$ is called Lie bialgebra. The relation between Poisson Lie groups and Lie bialgebras has been proved by Drinfeld \cite{Dr1}:

\begin{theorem}\label{thm: dr}
If $(G,\pi_G)$ is a Poisson Lie group, then the linearization of $\pi_G$ at $e$ defines a Lie algebra structure on $\mathfrak{g}^*$ such that $(\mathfrak{g},\mathfrak{g}^*)$ form a Lie bialgebra over $\mathfrak{g}$,
called the tangent Lie bialgebra to $(G,\pi_G)$. Conversely, if $G$ is connected and simply connected, then every Lie bialgebra $(\mathfrak{g},\mathfrak{g}^*)$ over $\mathfrak{g}$ defines a unique multiplicative Poisson
structure $\pi_G$ on $G$ such that $(\mathfrak{g},\mathfrak{g}^*)$ is the tangent Lie bialgebra to the Poisson Lie group $(G,\pi_G)$.
\end{theorem}

From this theorem it follows that there is a unique connected and simply connected Poisson Lie group $(G^*,\pi_{G\st})$, called the dual of $(G,\pi_G)$, associated to the Lie bialgebra $(\mathfrak{g}^*,\delta)$. If $G$ is connected and simply connected, then the dual of $G^*$ is $G$.

\begin{example}[$\mathfrak{g}=ax+b$]\label{ex: 1}

Consider the Lie algebra $\mathfrak{g}$ spanned by $X$ and $Y$ with commutator
\begin{equation}
[X,Y]=Y
\end{equation}
and cobracket given by
\begin{equation}
\delta(X)=0 \quad \delta(Y)= X\wedge Y.
\end{equation}
The Lie bracket on $\mathfrak{g}^*$ is given by
$$
[X^*,Y^*]=Y^*.
$$
A matrix representation of $\mathfrak{g}$ is the Lie algebra $\mathfrak{gl}(2,\mathbb{R})$ via
$$
X = \left(\begin{matrix} 1 & 0 \\ 0 & 0 \end{matrix}\right) \quad Y = \left(\begin{matrix} 0 & 1 \\ 0 & 0 \end{matrix}\right)
$$
and
$$
X^* = \left(\begin{matrix} 0 & 0 \\ 0 & 1 \end{matrix}\right) \quad Y^* = \left(\begin{matrix} 0 & 0 \\ 1 & 0 \end{matrix}\right)
$$
with the metric $\gamma(a,b)= tr(aJbJ)$ and $J=\left(\begin{matrix} 0 & 1 \\ 1 & 0 \end{matrix}\right)$.

The corresponding Poisson Lie group $G$ and dual $G^*$
are subgroups 
of $GL(2,\mathbb{R})$ of matrices with positive determinant are given by
\begin{equation}
G=\left\lbrace \left(\begin{matrix} 1 & 0 \\ \xi & \eta \end{matrix}\right)\; :\eta>0\right\rbrace \qquad G^*=\left\lbrace \left(\begin{matrix} s & t \\ 0 & 1 \end{matrix}\right)\; :s>0\right\rbrace
\end{equation}

\end{example}

\section{Poisson actions and Momentum maps}\label{sec_mm}

In this section we first introduce the concept of Poisson action of a Poisson Lie group on a Poisson manifold, which generalizes the canonical action of a Lie group on a symplectic manifold. We define momentum maps associated to such actions and finally we consider the particular case of a Poisson Lie group $G$ acting on its dual $G^*$ by dressing transformations. This allows us to study the symplectic leaves of $G$ that are exactly the orbits of the dressing action. These topics can be found e.g. in \cite{Lu3}, \cite{Lu1} and \cite{STS}.

%Recall that a canonical action of a Lie group $G$ on a Poisson manifold $M$ is defined as a group action which preserves the Poisson structure. 
%On the other side, a Poisson action is an action of a Poisson Lie group on a Poisson manifold satisfying a different property of compatibility between the Poisson bivectors of both manifolds. When the Poisson structure on the Lie group is trivial we recover the canonical actions.
 
%metti referenze qua

From now on we assume that $G$ is connected and simply connected.
\begin{definition}
 The action $\Phi:G\times M\rightarrow M$ of a Poisson Lie group $(G,\pi_G)$ on a Poisson manifold $(M,\pi)$ is called Poisson action if $\Phi$ is a Poisson map, where $G\times M$ is a Poisson  manifold with structure $\pi_G\oplus\pi$.
\end{definition}
This definition generalizes the notion of canonical action; indeed, if 
$G$ carries the trivial Poisson structure $\pi_G=0$, the action $\Phi$ is Poisson if and only if it preserves $\pi$, i.e. if it is canonical. In general, the structure $\pi$ is not invariant with respect to the action $\Phi$. The easiest examples of Poisson actions are given by the left and right actions of $G$ on itself.

For an action $\Phi: G\times M \to M$ we use $\Phi: \mathfrak{g} \to Vect\, M: X \mapsto \Phi_{X}$ to denote the Lie algebra anti-homomorphism which defines the infinitesimal generators of this action. The proof of the following Theorem can be found in \cite{LuWe1}.
\begin{theorem}
The action $\Phi: G\times M\rightarrow M$ is a Poisson action if and
only if
\begin{equation}\label{eq: pa}
L_{\Phi_{X}}(\pi)=(\Phi\wedge\Phi)\delta(X)
\end{equation}
 for any $X\in\mathfrak{g}$, where $L$ denotes the Lie derivative and $\delta$ is the derivative of $\pi_{G}$ at $e$.
\end{theorem}

Let $\Phi:G\times M \to M$ be a Poisson action of $(G, \pi_{G})$ on $(M,\pi)$. Let $G^*$ be the dual Poisson Lie group of $G$ and let $\Phi_{X}$ be the vector field on $M$ which generates the action $\Phi$.
In this formalism the definition of momentum map reads (Lu, \cite{Lu3}, \cite{Lu1}):

\begin{definition}\label{def: mm}
A momentum map for the Poisson action $\Phi:G\times M\rightarrow M$ is a map $\boldsymbol{\mu}: M\rightarrow G^*$ such that
\begin{equation}\label{eq: mmp}
\Phi_{X}=\pi^{\sharp}(\boldsymbol{\mu}^*(\theta_{X}))
\end{equation}
where $\theta_{X}$ is the left invariant 1-form on $G^*$ defined by the element $X\in\mathfrak{g}=(T_eG^*)^*$ and $\boldsymbol{\mu}^*$ is the cotangent lift $T^* G^*\rightarrow T^*M$.
\end{definition}
In other words, the momentum map generates the vector field $\Phi_{X}$ via the construction
$$
X\in \mathfrak{g} \to \theta_{X}\in T^*G^* \to \alpha_{X}= \boldsymbol{\mu}^*(\theta_{X})\in T^*M \to \pi^{\sharp}(\alpha_{X})\in TM
$$
It is important to remark that Noether's theorem still holds in this general context.
\begin{theorem}
Let $\Phi:G\times M \to M$ a Poisson action with momentum map $\boldsymbol{\mu}: M\rightarrow G^*$. If $H\in C^{\infty}(M)$ is $G$-invariant, then $\boldsymbol{\mu}$ is an integral of the Hamiltonian vector field associated to $H$.
\end{theorem}

It is important to point out that in this setting the vector field $\Phi_{X}$ is not Hamiltonian, unless the Poisson structure on $G$ is trivial. In this case $G^*=\mathfrak{g}^*$, the differential 1-form $\theta_{X}$ is the constant 1-form $X$ on $\mathfrak{g}^*$, and
\begin{equation}
\boldsymbol{\mu}^*(\theta_{X})=d(H_{X}),\quad\text{where}\quad H_{X}(m)=\langle\boldsymbol{\mu}(m),X \rangle.
\end{equation}
This implies that the momentum map is the canonical one and 
\begin{equation}
\Phi_{X}=\pi^{\sharp}(dH_{X})=\{H_{X}, \cdot\}.
\end{equation}
In other words, $\Phi_{X}$ is the Hamiltonian vector field with Hamiltonian $H_{X}\in C^{\infty}(M)$. We observe that,
when $\pi_G$ is not trivial, $\theta_{X}$ is a Maurer-Cartan form, hence $\boldsymbol{\mu}^*(\theta_{X})$ can not be written as a differential of a Hamiltonian function.
In the following we give an example for the infinitesimal generator in this general case.

\subsection{Dressing Transformations}\label{sec: dressing}

One of the most important example of Poisson action is the dressing action of $G$ on $G^*$. The name ``dressing'' comes from the theory of integrable systems and was introduced in this context in \cite{STS}. Interesting examples can be found in \cite{AM}.
We remark that, given a Poisson Lie group $(G, \pi_{G})$, the left (right) invariant 1-forms on $G^*$ form a Lie algebra with respect to the bracket:
$$
[\alpha,\beta] = L_{\pi^{\sharp}(\alpha)}\beta-L_{\pi^{\sharp}(\beta)}\alpha - d(\pi(\alpha, \beta)). 
$$

For $X\in\mathfrak{g}$, let $\theta_{X}$ be the left invariant 1-form on $G^*$ with value $X$ at $e$. Let us define the vector field on $G^*$

\begin{equation}\label{eq: idr}
	l(X)=\pi_{G\st}\sh(\theta_{X}).
\end{equation}
The map $l: \mathfrak{g}\to TG^*: X\mapsto l(X)$ is a Lie algebra anti-homomorphism. We call $l$ the left infinitesimal dressing action of $\mathfrak{g}$ on $G^*$; its linearization at $e$ is the coadjoint action of $\mathfrak{g}$ on $\mathfrak{g}^*$. Similarly we can define the right infinitesimal dressing action.

Let $l(X)$ (resp. $r(X)$) a left (resp. right) dressing vector field on $G^*$. If all the dressing vector fields are complete, we can integrate the $\mathfrak{g}$-action into an action of $G$ on $G^*$ called the
dressing action and we say that the dressing actions consist of dressing transformations. 
\begin{definition}
A multiplicative Poisson tensor $\pi_G$ on $G$ is complete if each left (equiv. right) dressing vector field is complete on $G$.
\end{definition}

From the definition of dressing action follows (the proof can be found in \cite{STS}) that the orbits of the right or left dressing action of $G^*$ (resp. $G$) are 
the symplectic leaves of $G$ (resp. $G^*$).

It can be proved (see \cite{Lu3}) that if $\pi_{G}$ is complete, both left and right dressing actions are Poisson actions with momentum map given by the identity.

Assume that $G$ is a complete Poisson Lie group. We denote respectively the left (resp. right) dressing action of $G$ on its dual $G^*$ by $g\mapsto l_g$ (resp. $g\mapsto r_g$).

\begin{definition}
A momentum map $\boldsymbol{\mu}:M\rightarrow G^*$ for a left (resp. right) Poisson action $\Phi$ is called G-equivariant if it is such with respect to the left dressing action of $G$ on $G^*$, that is, 
$\boldsymbol{\mu}\circ \Phi_g=\lambda_g\circ \boldsymbol{\mu}$ (resp. $\boldsymbol{\mu}\circ \Phi_g=\rho_g\circ \boldsymbol{\mu}$)
\end{definition}
It is important to remark that a momentum map is $G$-equivariant if and only if it is a Poisson map, i.e. $\boldsymbol{\mu}_*\pi=\pi_{G^*}$.
\begin{definition}
 An action $\Phi: G\times M \to M$ of a Poisson Lie group $(G, \pi_{G})$ on a Poisson manifold $(M, \pi)$ is said Hamiltonian if it is a Poisson action generated by an equivariant momentum map.
\end{definition}
%esempio?

\section{Poisson Reduction}\label{sec: pr}
\label{sec: poisson reduction}

In this section we present the main result of this paper. We show that, given a Hamiltonian action $\Phi$, as defined above, we can define a reduced manifold in terms of momentum map and prove that it is a Poisson manifold. The approach used is a generalization of the orbit reduction \cite{Ml} in symplectic geometry. Recall that, under certain conditions, the orbit space of $\Phi$ is a smooth manifold and it carries a Poisson structure. First, we give an alternate proof of this claim. Then, we consider a generic orbit $\mathcal{O}_{u}$ of the dressing action of $G$ on $G^*$, for $u\in G^*$, and we prove that the set $\boldsymbol{\mu}^{-1}(\mathcal{O}_{u})/G$ is a regular quotient manifold with Poisson structure induced by the Poisson structure on $M$. Similarly to the symplectic case, this reduced space is isomorphic to the space $\boldsymbol{\mu}^{-1}(u)/G_{u}$ which will be regarded as the Poisson reduced space.

\subsection{Poisson structure on $M/G$}

Consider a Hamiltonian action of a connected and simply connected Poisson Lie group 
$(G,\pi_{G})$ on a Poisson manifold $(M,\pi)$. It is known that, if the action is proper and free, the orbit space $M/G$ is a smooth manifold, it carries a Poisson structure such that the natural projection $M \to M/G$ is a Poisson map (a proof of this result can be found in \cite{STS}).
In this section we give an alternate proof of this result, by introducing an explicit formulation for the infinitesimal generator of the Hamiltonian action, in terms of local coordinates. 

As discussed in the previous section, a Hamiltonian action is a Poisson action induced by an equivariant momentum map $\boldsymbol{\mu}: M \to G^*$ by formula (\ref{eq: mmp}). In other words, the map
$$
\alpha: \mathfrak{g} \to \Omega^1(M): X \mapsto \alpha_{X}=\boldsymbol{\mu}^*(\theta_{X})
$$
is a Lie algebra homomorphism such that
$$
\Phi_{X}=\pi^{\sharp}(\alpha_{X})
$$
The dual map of $\alpha$ defines a $\mathfrak{g}^*$-valued 1-form on $M$, still denoted by $\alpha$, satisfying Maurer-Cartan equation (as proved in \cite{Lu3})
$$
d\alpha+\frac{1}{2}[\alpha,\alpha]_{\mathfrak{g}^*}=0.
$$
In particular,
$$
\{\alpha_{X}: X\in\mathfrak{g}\}
$$
defines a foliation $\mathcal{F}$ on $M$.

\begin{lemma}\label{thm: m/g}
The space of $G$-invariant functions on $M$ is closed under Poisson bracket. Hence $\pi$ defines a Poisson structure on $M/G$
\end{lemma}

\begin{proof}
Let $H_{i}$, $i=1, \dots n$ be local coordinates on $M$ such that 
$$
\mathcal{F}= Ker\{dH_{1}, \dots, dH_{n}\}
$$
Then 
\begin{equation}\label{eq: al}
\alpha_{X}=\sum_{i}c_{i}(X)dH_i
\end{equation}
and
\begin{equation}\label{eq: xis}
\Phi_{X}[f]=\pi^{\sharp}(\alpha_{X})=\sum_{i}c_{i}(X)\lbrace H_j,f\rbrace_{M}.
\end{equation}
This implies that a function $f\in C^{\infty}(M)$ is $G$-invariant ($\Phi_{X}[f]=0$) if and only if 
$\lbrace H_i,f\rbrace=0$ for any $i$. If $f,g$ are $G$-invariant functions on $M$, we have $\lbrace H_i,f\rbrace=\lbrace H_i,g\rbrace=0$ for any $i$. Then, using the Jacobi identity we get
 $\lbrace H_i,\lbrace f,g\rbrace\rbrace=0$. Since $G$ is connected, the result follows.
\end{proof}

\subsection{Poisson reduced space}

Assume that $G$ is connected, simply connected and complete. In order to define a reduced space and to prove that it is a Poisson manifold we consider a generic orbit $\mathcal{O}_u$ of the dressing orbit of $G$ on $G^*$ passing through $u\in G^*$. First, we prove the following:

\begin{lemma}
Let $\Phi:G\times M \to M$ be a free and Hamiltonian action of a compact Poisson Lie group $(G,\pi_{G})$ on a Poisson manifold $(M,\pi)$. Then:
\begin{itemize}
	\item[(i)] $\mathcal{O}_u$ is closed and the Poisson structure $\pi_{G^*}$ does not depend on the transversal coordinates on $\mathcal{O}_u$.
	\item[(ii)] $\boldsymbol{\mu}^{-1}(\mathcal{O}_u)/G$ is a smooth manifold.
\end{itemize}
\end{lemma}

\begin{proof}

\begin{itemize}
	\item[(i)] If $G$ is compact, any $G$-action is automatically proper. This implies that, given $u\in G^*$ the generic orbit $\mathcal{O}_u$ of the dressing action is closed. From section (\ref{sec: dressing}) we know that $\mathcal{O}_u$ is the symplectic leaf through $u$. Using the local description of Poisson manifolds introduced in Theorem (\ref{thm: split}) it is evident that $\pi_{G^*}$ restricted to $\mathcal{O}_u$ does not depend on the transversal coordinates $y_{i}$.
	\item[(ii)] If the action $\Phi$ is free, the momentum map $\boldsymbol{\mu}:M \to G^*$ is a submersion onto some open subset of $G^*$. This implies that $\boldsymbol{\mu}^{-1}(u)$ is a closed submanifold of $M$. As $\boldsymbol{\mu}$ is equivariant, it follows that $\boldsymbol{\mu}^{-1}(u)$ is $G$-invariant. Free and proper actions of $G$ on $M$ restrict to free and proper $G$-actions on $G$-invariant submanifolds. In particular, the action of $G$ on $\boldsymbol{\mu}^{-1}(u)$ is still proper, then $G\cdot \boldsymbol{\mu}^{-1}(u)$ is closed. Using the equivariance we have that $G\cdot \boldsymbol{\mu}^{-1}(u)= \boldsymbol{\mu}^{-1}(\mathcal{O}_u)$, which is still $G$-invariant. The action of $G$ on $\boldsymbol{\mu}^{-1}(\mathcal{O}_u)$ is proper and free, so we can conclude that the orbit space $\boldsymbol{\mu}^{-1}(\mathcal{O}_u)/G$ is a smooth manifold.
\end{itemize}

\end{proof}

We aim to prove that the manifold $N/G:=\boldsymbol{\mu}^{-1}(\mathcal{O}_u)/G$ carries a Poisson structure. 
In the previous Lemma we stated that $\pi_{G^{*}}$ restricted to $\mathcal{O}_u$ does not depend on the transversal coordinates $y_i$'s; if $x_{i}$ are local coordinates along $N=\boldsymbol{\mu}^{-1}(\mathcal{O}_u)$ and $H_{i}$ are pullback of the transversal coordinates $y_{i}$'s by 
\begin{equation}
H_{i}:= y_{i}\circ \boldsymbol{\mu}
\end{equation}
we can easily deduce that the Poisson structure $\pi$ on $M$ involves derivatives in $H_{i}$ only in the combination
$$
\partial_{x_i}\wedge\partial_{H_i}
$$
This is evident because the differential $d\boldsymbol{\mu}$ between $TM\vert_N/TN$ and $TG^{*}/T\mathcal{O}_u$ is a bijective map. Moreover, since $\{y_{i},y_{j}\}$ vanishes on the orbit $\mathcal{O}_u$, $\{H_{i},H_{j}\}$ vanishes on the preimage $N$ and $dH_{i}$'s are in the span of $\{\alpha_{X}:X\in\mathfrak{g}\}$.

Now we introduce the ideal $\mathcal{I}$ generated by $H_i$ and prove some properties.
\begin{lemma}
Let $\mathcal{I}=\{f\in C^{\infty}(M): f\vert_{N}=0\}$. 
\begin{itemize}
	\item[(i)] $\mathcal{I}$ is defined in an open $G$-invariant neighborhood $U$ of $N$.
	\item [(ii)] $\mathcal{I}$ is closed under Poisson bracket.
\end{itemize}
\end{lemma}

\begin{proof}
\begin{itemize}
	\item[(i)] The coordinates $H_i$ are locally defined but we can show that $\mathcal{I}$ is globally defined.
Considering a different neighborhood on the orbit of $G^{*}$ we have transversal coordinates $y_i^{\prime}$ and their pullback to $M$ will be $H_i^{\prime}=y_i^{\prime}\circ \boldsymbol{\mu}$.
The coordinates $H_i^{\prime}$ are defined in a different open neighborhood $V$ of $N$, but we can see that the ideal $\mathcal{I}$ generated by $H_i$ coincides with $\mathcal{I}^{\prime}$ generated by $H_i^{\prime}$ 
on the intersection of $U$ and $V$, then it is globally defined.
	\item [(ii)] Since $\boldsymbol{\mu}$ is a Poisson map we have:
$$
 \{ H_i,H_j\}_{M}=\{ y_i\circ \boldsymbol{\mu},y_j\circ \boldsymbol{\mu}\}_{M}=\{ y_i,y_j\}_{G^*}\circ \boldsymbol{\mu}.
$$
Hence the ideal $\mathcal{I}$ is closed under Poisson brackets.
\end{itemize}
\end{proof}

Motivated by this Lemma we use the following identification
$$
C^{\infty}(N/G)\simeq(C^{\infty}(U)/\mathcal{I})^G.
$$

\begin{lemma}\label{lem: id2}
 Suppose that $N/G$ is an embedded submanifold of the smooth manifold $M/G$, then
\begin{equation}
 (C^{\infty}(U)/\mathcal{I})^G \simeq(C^{\infty}(U)^G + \mathcal{I})/\mathcal{I}
\end{equation}
\end{lemma}

\begin{proof}
 Let $f$ be a smooth function on $U$ satisfying $[f]\in (C^{\infty}(U)/\mathcal{I})^G$. As the equivalence class $[f]$ is $G$-invariant, we have
\begin{equation}
f(G\cdot m)=f(m)+i(m),
\end{equation}
where $i\in \mathcal{I}$ and $G\cdot m$ is a generic orbit of the Hamiltonian action of $G$ on $M$. It is clear that $f\vert_N$ is $G$-invariant and hence it defines a smooth function $\bar{f}\in C^{\infty}(N/G)$.
Since $N/G$ is a $k$-dimensional embedded submanifold of the $n$-dimensional smooth manifold $M/G$, the inclusion map $\iota: N/G\rightarrow M/G$ has local coordinates representation:
\begin{equation}
 (x_1,\dots,x_k)\mapsto (x_1,\dots,x_k,c_{k+1},\dots,c_n)
\end{equation}
where $c_i$ are constants. Hence we can extend $\bar{f}$ to a smooth function $\phi$ on $M/G$ by setting $\bar{f}(x_1,\dots,x_k)=\phi(x_1,\dots,x_k,0,\dots,0)$.
The pullback $\tilde{f}$ of $\phi$ by $\text{pr}:M\rightarrow M/G$ is $G$-invariant and satisfies
\begin{equation}
 \tilde{f}-f\vert_N=0,
\end{equation}
hence $\tilde{f}-f\in \mathcal{I}$.
\end{proof}

Using these results we can prove the following:

\begin{theorem}\label{thm: pred}
Let $\Phi:G\times M\rightarrow M$ be a free Hamiltonian action of a compact Poisson Lie group $(G,\pi_G)$ on a Poisson manifold $(M,\pi)$ with momentum map $\boldsymbol{\mu}:M\rightarrow G^*$. The orbit space $N/G$ has a Poisson structure induced by $\pi$.
\end{theorem}

\begin{proof}
First we prove that the Poisson bracket of $M$ induces a well defined Poisson bracket on $(C^{\infty}(U)^G+\mathcal{I})/\mathcal{I}$.
In fact, for any $f+i\in C^{\infty}(U)^{G}/\mathcal{I}$ and $j\in \mathcal{I}$ the Poisson bracket $\{ f+i,j\}$ still belongs to the ideal $\mathcal{I}$.
Since the ideal $\mathcal{I}$ is closed under Poisson brackets, $\{ i,j\}$ belongs to $\mathcal{I}$.
The function $j$, by definition on the ideal $\mathcal{I}$, can be written as a linear combination of $H_i$, so $\{ f,j\}=\sum_i a_i\{ f,H_i\}$.
By Lemma \ref{thm: m/g}, we have $\{ f,H_i\}=0$, hence $\{ f+i,j\}\in \mathcal{I}$ as stated.
Finally, using the isomorphism proved in the Lemma (\ref{lem: id2}) and the identification $
C^{\infty}(N/G)\simeq(C^{\infty}(U)/\mathcal{I})^G$, the claim is proved.
\end{proof}

Finally, we observe that there is a natural isomorphism
\begin{equation}
 \boldsymbol{\mu}^{-1}(u)/G_{u}\simeq \boldsymbol{\mu}^{-1}(\mathcal{O}_u)/G.
\end{equation}
We refer to $\boldsymbol{\mu}^{-1}(u)/G_{u}$ as the Poisson reduced space.

\section{An example}
\label{sec: ex}
In this section we discuss a concrete example of Poisson reduction. Consider the Lie bialgebra $\mathfrak{g}=ax+b$ discussed in Example (\ref{ex: 1}). The Poisson tensor on the dual Poisson Lie group $G^*$ is given, in the coordinates $(s,t)$ introduced in the matrix representation, by
\begin{equation}
\pi_{G^*}=st\partial_s\wedge\partial_t.
\end{equation}
It is clear that $(s,t)$ are global coordinates on $G^*$. 
First, we need to study the orbits of the dressing action.
Remember that the dressing orbits $\mathcal{O}_u$ through a point $u\in G^*$ are the same as the symplectic leaves, hence it 
is clear that they are determined by the equation $t=0$. The symplectic foliation of the manifold $G^*$ in this case is given by two open orbits, determined by the conditions
$t>0$ and $t<0$ respectively, and a closed orbit given by $t=0$ and $a\in\mathbb{R}^+$.

Consider a Hamiltonian action $\Phi:G\times M \to M$ of $G$ on a generic Poisson manifold $M$ induced by the equivariant momentum map $\boldsymbol{\mu}:M\rightarrow G^*$.
Its pullback
\begin{equation}
\boldsymbol{\mu}\st: C^{\infty}(G^*)\longrightarrow C^{\infty}(M)
\end{equation}
maps the coordinates $s$ and $t$ on $G^*$ to 
$$
x(u)=s(\boldsymbol{\mu}(u)) \qquad y(u)=t(\boldsymbol{\mu}(u)).
$$  
It is important to underline that we have no information on the dimension of $M$, so $x$ and $y$ are just a pair of the possible coordinates. Nevertheless, since $\boldsymbol{\mu}$ is a Poisson map,  we have
\begin{equation}
\lbrace x,y\rbrace=xy
\end{equation}
on $M$. The infinitesimal generators of the action $\Phi$ can be written in terms of these 
coordinates $(x,y)$ as
\begin{equation}
\Phi(X)=x\{y,\cdot \}\quad \Phi(Y)=x\{x^{-1},\cdot \}.
\end{equation}

In the following, we discuss the Poisson reduction case by case, by considering the different dressing orbits studied above.

\paragraph{Case 1: $(t>0)$} Consider the dressing orbit $\mathcal{O}_u$ determined by the condition $t>0$. Since $s$ and $t$ are both positive, we can put
\begin{equation}
  x=e^p,\quad y=e^q.
\end{equation}
Since $\lbrace x,y\rbrace=xy$ we have
\begin{equation}
	\lbrace p,q\rbrace=1.
\end{equation} 
For this reason the preimage of the dressing orbit can be split as $N=\mathbb{R}^2\times M_1$ and $C^{\infty}(N)$ is given explicitly by the set of functions generated by $y^{-1}$. The infinitesimal generators are given by
\begin{equation}
\Phi(X)=e^p\{e^q,\cdot \}\qquad \Phi(Y)=e^p\{e^{-p},\cdot \}
\end{equation}
which is the action of $G$ on the plane. Hence the Poisson reduction in this case is given by
\begin{equation}
 (C^{\infty}(M)[y^{-1}])^G.
\end{equation}

\paragraph{Case 2: $(t<0)$} This case is similar, with the only difference that $y=-e^q$. 

\paragraph{Case 3: $(t=0)$} The orbit $\mathcal{O}_u$ is given by fixed points on the line $t=0$, then we choose the point $s=1$. Consider the ideal $\mathcal{I}=\langle x-1,y\rangle$ of functions vanishing on $N$. It is easy to check that it is $G$-invariant, hence the Poisson reduction in this case is simply given by
\begin{equation}
 (C^{\infty}(M)/\mathcal{I})^G.
\end{equation}

\subsection{Questions and future directions}

The theory of Poisson reduction can be further developed, as it has been obtained under the assumption that the orbit space $M/G$ is a smooth manifold. This result could be proved under weaker hypothesis, for instance requiring that $M/G$ is an orbifold. 

As stated in the introduction, the idea of momentum map and Poisson reduction can be also used for the study of symmetries in quantum mechanics. In particular, the approach of deformation quantization would provide a relation between classical and quantum symmetries. A notion of quantum momentum map has been defined in \cite{me}, \cite{me1} and it can be used to define the quantization of the Poisson reduction.

At classical level, Poisson reduction could be generalized to actions of Dirac Lie groups \cite{MJ} on Dirac manifolds \cite{Co}. Finally, a possible development of this theory is its integration to symplectic groupoids by means of the theories on the integrability of Poisson brackets \cite{CF} and Poisson Lie group actions \cite{FP}.

%% The Appendices part is started with the command \appendix;
%% appendix sections are then done as normal sections
%% \appendix

%% \section{}
%% \label{}

%% References
%%
%% Following citation commands can be used in the body text:
%% Usage of \cite is as follows:
%%   \cite{key}          ==>>  [#]
%%   \cite[chap. 2]{key} ==>>  [#, chap. 2]
%%   \citet{key}         ==>>  Author [#]

%% References with bibTeX database:

\bibliographystyle{alpha}

\end{document}